\documentclass[11pt]{article}
\usepackage[dvips]{epsfig}
\usepackage{color}

\newcommand{\F}{\mathcal{F}}    
\newcommand{\Pa}{\mathrm{Pa}}  
\newcommand{\CA}{\mathcal{A}}

\newcommand{\DC}{\mathcal{D}} 
\newcommand{\CH}{\mathcal{H}}
\newcommand{\TO}{\mathrm{T0}}   
\newcommand{\TS}{\mathrm{TS}}   
\newcommand{\TSO}{\mathrm{TS0}}
\newcommand{\CS}{\mathrm{CS}}
\newcommand{\HS}{\mathrm{HS}}
\def\Fix{\mathrm{Fix}}
\def\fix{\mathrm{fix}}
\def\da{\mathrm{da}}
\def\v{\mathrm{v}}
\def\ds{\mathrm{ds}}
\def\h{\mathrm{h}}
\def\r{\mathrm{r}}
\def\id{\mathrm{id}}
\def\D6{\mathcal{D}_6}
\def\Dn{\mathcal{D}_n}
\def\P{\mathcal{P}}
\def\Q{\mathcal{Q}}
\def\R{\mathcal{R}}

\title{Enumeration of Symmetry Classes of Convex Polyominoes on the Honeycomb Lattice
\footnote {With the partial support of CNRS (France), NSERC (Canada) and FCAR (Qu\'ebec).
This is the full version of a paper presented at the FPSAC Conference in Vancouver, Canada,
June 28 -- July 2, 2004 (see \cite{resumeFPSAC04}).}}
\author{Dominique Gouyou-Beauchamps, LRI, CNRS, France\\ and Pierre Leroux, LaCIM, UQAM, Canada}
\date{March 9, 2004}   
\begin{document}         
\maketitle                 
%
%
%
\begin{abstract} 
\emph{Hexagonal} polyominoes are polyominoes on the honeycomb lattice.
We enumerate the symmetry classes of \emph{convex} hexagonal polyominoes.
Here \emph{convexity} is to be understood as convexity along the three
main column directions.
We deduce the generating series of \emph{free} (i.e. up to reflection and rotation) 
and of \emph{asymmetric} convex hexagonal polyominoes, according to area and half-perimeter.
We give explicit formulas or implicit functional equations for the generating series,
which are convenient for computer algebra.
\end{abstract}
%
%
\section{Introduction}

\begin{figure}[h]
\epsfysize=5.0cm  
\centerline{\epsfbox{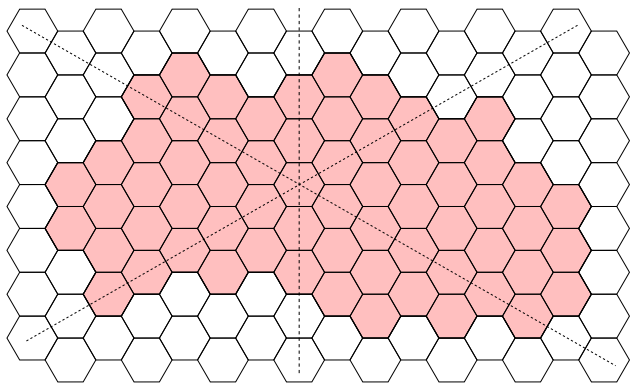}}
\vspace{-0.45\baselineskip}
\caption{A convex polyomino and its convexity directions}
\label{fig:convex}
\end{figure}

An \emph{hexagonal polyomino} is a finite connected set of basic cells of the
honeycomb lattice in the plane.	 
Note that the hexagons of our lattice have two sides parallel to the horizontal axis.
See Figure \ref{fig:convex}.
Unless otherwise stated, all the polyominoes considered here are hexagonal.
The \emph{area} of a polyomino is the number of the cells composing it.
Its \emph{perimeter} is the number of line segments on its boundary.
We say that a polyomino is \emph{convex along a direction} if 
the intersection with any line parallel to this direction and passing through
the center of a cell is connected.
The convexity directions are characterized by the angle $\alpha$ ($0\leq\alpha\leq\pi$)
which they form with the positive horizontal axis.

Various convexity concepts have been introduced in the literature 
for hexagonal polyominoes, depending on the required convexity directions.
Following the nomenclature of Denise, D\"ur, and Hassani  \cite{ddi97},
we mention the \emph{$EG$-convex} polyominoes, where $\alpha = 0$ and $\pi/2$, studied
by Guttmann and Enting \cite{ge88} and by Lin and Chang \cite{lc88},
the \emph{$C^1$-convex} polyominoes, where $\alpha = \pi/2$, enumerated according to many
parameters by Lin and Wu \cite{lw90} and by Fereti\'c and Svrtan \cite {fs93},
the \emph{strongly convex} polyominoes, 
where $\alpha = 0$, $\pi/3$ and $2\pi/3$, 
introduced by Hassani \cite{ibn96} and studied in \cite{ibn96} and \cite{ddi97},
and finally the $C$- or \emph{$C^3$-convex} polyominoes, 
where $\alpha= \pi/6$, $\pi/2$ and $5\pi/6$,
introduced and enumerated according to perimeter in \cite{ibn96} and \cite{ddi97}.
In particular, Hassani gives explicitly the algebraic generating function for 
$C$-convex polyominoes according to half perimeter.

It is this last class that interest us here, and that we  call
simply \emph{convex polyminoes}. See Figure \ref{fig:convex} for an example.
This concept is a natural extension of (row and column) convexity on the
square lattice. 

These polyominoes are traditionally taken up to translation.
However it is natural to consider them also up to rotation and reflection,
as objects living freely in space. Following V\"oge, Guttmann and Jensen \cite{vgj02}, 
we call these equivalence classes \emph{free polyominoes}.
In organic chemistry, free polyminoes represent benzenoid hydrocarbons.
See \cite{vgj02} where these molecules (without the convexity
property) are enumerated by an exhaustive generation approach. 

Our objective is to enumerate free convex polyominoes, according to area and half perimeter.
Following the approach of Leroux, Rassart and Robitaille \cite{lrr98} 
for the square lattice,
we consider them as orbits of the dihedral group $\D6$ 
(the group of isometries of the regular hexagon), acting on convex polyominoes.
The Cauchy-Frobenius Formula (alias Burnside's Lemma) can be used to count the orbits and
we are thus lead to enumerate the symmetry classes 
$\Fix(h)$ of convex polyominoes, for each element $h$ of the group $\D6$, 

It is also possible to enumerate convex polyominoes which are asymmetric or which have
exactly the symmetries of a given subgroup $H$ of $\D6$, using M\"obius inversion
in the lattice of subgroups of $\D6$. For this purpose, we  also enumerate
the symmetry classes $F_{\geq H}$ for each subgroup $H$ of $\D6$.

For any class $\F$ of convex polyominoes, we  denote its generating series 
by $\F(x,q,u,v,t)$, where the variable $x$ marks the number of columns, $q$ marks the area,
$u$ marks the size of the  first column (on the left), $v$, the  size of the last column,
and $t$ the half perimeter; for example, the polyomino of Figure \ref{fig:convex}
has weight $x^{14}q^{64}u^2v^3t^{35}$.
It is possible that some variables do not appear in some generating series. 
The generating series will be given by explicit formulas 
or implicit functional equations, which are convenient for computer algebra.
%
\section{Preliminaries}
%
\subsection{Particular classes of convex polyominoes}
Some familiar classes of convex polyominoes of the square lattice are naturally
found on the honeycomb lattice and are useful. 
It is the case notably of \emph{partition} and \emph{staircase} (or \emph{parallelogram}) 
polyominoes and of a variant of \emph{stack} polyominoes.
\begin{figure}[h]
\epsfysize=4cm  
\centerline{\epsfbox{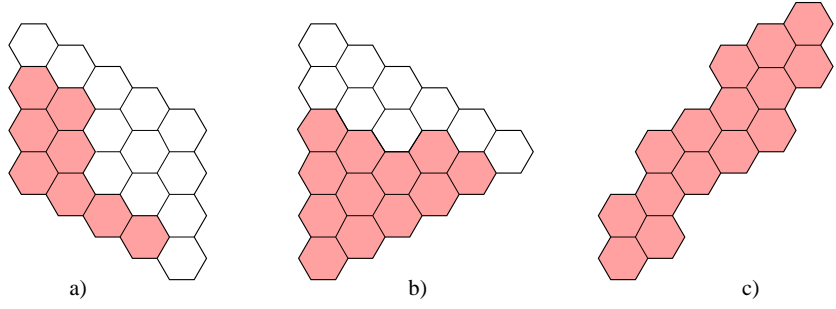}}
\vspace{-0.45\baselineskip}
\caption{Partition and staircase polyominoes
\label{fig:partapara}}
\end{figure}
%
\subsubsection{Partition polyominoes}
Figure \ref{fig:partapara}a represents the partition $(4,2,2)$ contained in a 
\emph{rectangle} of size $5\times4$ in the honeycomb lattice. 
Figure \ref{fig:partapara}b represents the distinct part partition $(5,4,2,1)$,
with parts bounded by 6. 
We denote by $D_m(u,q)$ the generating polynomial of distinct part partitions 
with parts bounded by $m$. 
Here the  variable $u$ marks the number of parts. We have
\begin{equation}
D_m(u,q)=(1+uq)(1+uq^2)\cdots (1+uq^m)\mbox{~~and~~}D_0(u,q)=1.
\label{eq:partadist}
\end{equation}
%
\subsubsection{Staircase polyominoes}
Figure \ref{fig:partapara}c represents a staircase polyomino 
from the square lattice (see for example \cite{bou93} or \cite{lr01})
redrawn on the honeycomb lattice. Observe that the half perimeter
is equal to $2p-1$ where $p$ is the half perimeter on the square lattice.
We know that these polyominoes are enumerated according to half perimeter by the 
Catalan numbers and according to area by the sequence M1175 of \cite{sp95} (A006958 of \cite{sl03})
whose generating series is a quotient of two $q$-Bessel functions.
 
We denote by $\Pa$, the set of staircase polyominoes (Pa for \emph{Parallelogram})
on the honeycomb lattice and by $\Pa(x,q,u,v,t)$, their generating series.
An analysis of the situation where a column is added on the right, following the 
method of M. Bousquet-Melou (compare with \cite{bou93}, Lemma 3.1), gives, for $\Pa(v)=\Pa(x,q,u,v,t)$,
\begin{equation}
\Pa(v)=\frac{xquvt^3}{1-quvt^2}+\frac{xqvt^2}{(1-qvt^2)(1-qv)}
\left( \Pa(1)-\Pa(vq)\right)
\label{eq:parafonct}
\end{equation}
so that
\begin{equation}
\Pa(v)=
\frac{J_1(1)+J_1(v)J_0(1)-J_1(1)J_0(v)}{J_0(1)},
\label{eq:byaformula}
\end{equation}
where
$$ J_1(v) 
=\sum_{n\geq 0}(-1)^n\frac{x^{n+1}v^{n+1}ut^{2n+3}q^{{n+2}\choose 2}}
{(qvt^2;q)_n(qv;q)_n(1-q^{n+1}uvt^2)}$$
and
$$J_0(v)=\sum_{n\geq 0}(-1)^{n}\frac{x^{n}v^nt^{2n}q^{{n+1}\choose 2}}
{(qvt^2;q)_{n}(qv;q)_{n}}.$$
Here we have used the familiar notation $(a;q)_n=(1-a)(1-aq)\cdots(1-aq^{n-1}$. We set
%
\begin{equation}
\Pa(x,q,u,v,t)=\sum_{i\geq 1\,j\geq 1}\Pa_{i,j}(x,q,t)u^iv^j.
\label{eq:Gij}
\end{equation}
%
\subsubsection{Stack polyominoes\label{sec:tas}}
\begin{figure}[h]
\epsfysize=4.0cm  
\centerline{\epsfbox{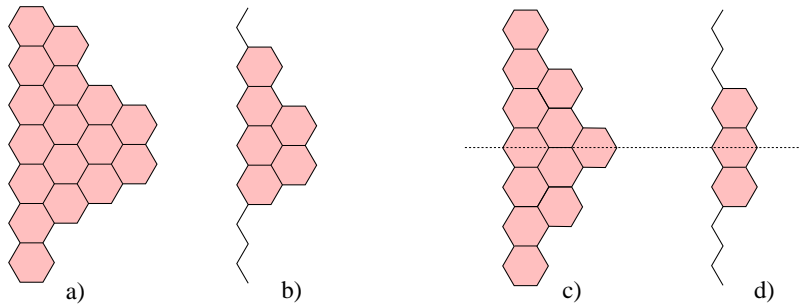}}
\vspace{-0.45\baselineskip}
\caption{Stacks and symmetric stacks
\label{fig:tas}}
\end{figure}

There exists a specific variant of stack polyominoes on the honeycomb lattice.
They consist of pyramidal stackings of hexagons, viewed sideways for our purposes. 
A first class (see Figure \ref{fig:tas}a), denoted by $T$ (for French \emph{tas}), appears 
in the literature under the name of \emph{pyramidal stacking of circles}; see \cite{prsv89}. 
Their generating series according to area is referenced as number M0687 in \cite{sp95} and 
A001524 in \cite{sl03}.

Let $T(x,u,q)$ be the generating series of stack polyominoes according to the number of columns
(the \emph{width}, marked by $x$), the size of the first column 
(the \emph{height}, marked by $u$) 
and the area, and let $T_n(x,q)=[u^n]T(x,u,q)$ be the  generating series of stacks whose
first column is of size $n$. Note that the half perimeter is equal to
twice the height plus the width so that the series $T(xt,ut^2,q)$ also 
keeps track of this parameter.

We have  
\begin{equation}
 T(x,u,q)=\sum_{m\geq 1}\frac{x^mq^{{m+1}\choose 2}u^m}
{\left((uq;q)_{m-1}\right) ^2(1-uq^m)}
\label{eq:tas}
\end{equation}
and
\begin{equation}
T_n(x,q)=\sum_{m=1}^nx^mq^{n+{m\choose 2}}\sum_{j=0}^{n-m}
\left[ _{~~m-1}^{m+j-1} \right]_q \left[ _{~m-2}^{n-j-2} \right]_q.
\label{eq:tasn}
\end{equation}
The polynomials $T_n(x,q)$ can also be rapidly computed by recurrence 
using the class $\TO_n$ of stack polyominoes whose first column is 
of size $n$, including empty cells at the two extremities. 
See Figure \ref{fig:tas}b.
Indeed, we  have
\begin{equation}
T_n(x,q)=xq^n\TO_{n-1}(x,q).
\label{eq:tasntaso}
\end{equation}
with $\TO_{0}(x,q)=1$, $\TO_{1}(x,q)=1+xq$, and, arguing on the existence of empty cells
at each extremity,
\begin{equation}
\TO_{n}(x,q)= (2+xq^n)\TO_{n-1}-\TO_{n-2},
\label{eq:recurrencandason}
\end{equation}
%
\subsubsection{Symmetric stacks}
Horizontally symmetric stacks (see Figures \ref{fig:tas}c 
and \ref{fig:tas}d), constitute the families $\TS$ and $\TSO$. 
Using the same notation as for stacks, we  have 
\begin{equation}
\TS(x,u,q) =\sum_{m\geq 1}\frac{x^mu^mq^{m(m+1)/2}(1+uq^m)}{(1-u^2q^2)(1-u^2q^4)\cdots
(1-u^2q^{2m})}.
\label{eq:tassym}
\end{equation}
Moreover, 
\begin{equation}
\TS_n(x,q)=xq^n\TSO_{n-1}(x,q).
\label{eq:tassymntassymo}
\end{equation}
where $\TSO_{0}(x,q) = \TSO_{-1}(x,q) = 1$, and
\begin{equation}
\TSO_{n}(x,q) = xq^n\TSO_{n-1}(x,q)+\TSO_{n-2}(x,q),
\label{eq:recurrencandassymon}
\end{equation}
%
\subsection{The dihedral group $\D6$}
\begin{figure}[h]
\epsfysize=4cm
\centerline{\epsfbox{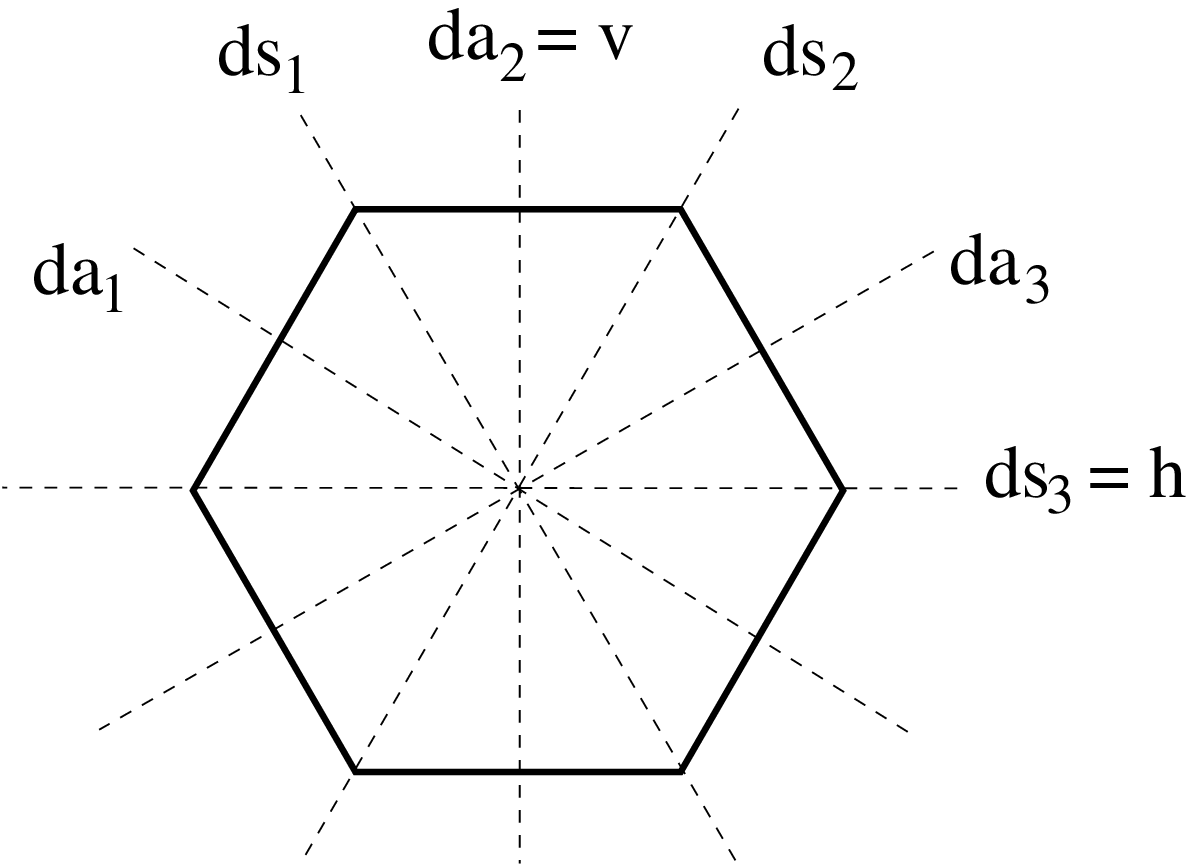}}
\vspace{-0.45\baselineskip}
\caption{The reflections of $\D6$}
\label{fig:axesdereflexion}
\end{figure}

The dihedral group $\D6$ is defined algebraically by 
$$\D6 = \langle\rho,\tau ~|~ \rho^6=1,~\tau^2=1,~\tau\rho\tau=\rho^{-1}\rangle.$$
Here $\D6$ is realized as the group of isometries of a regular hexagon,
with $\rho=\r$ = the (clockwise) rotation of $\pi/3$ radian and $\tau=\ds_3=\h$, 
the horizontal reflection. We have 
$$\D6=\{\id,\r,\r^2, \r^3, \r^4, \r^5,\da_1,\da_2,\da_3,\ds_1,\ds_2,\ds_3\},$$
where $\ds_2=\tau\r^2$, $\ds_1=\tau\r^4$, reflections according to vertex-vertex axes,
and $\da_3=\tau\r$, $\da_2=\tau\r^3$, and $\da_1=\tau\r^5$,
reflections according to the edge-edge axes.
See Figure \ref{fig:axesdereflexion}.

The dihedral group $\D6$ acts naturaly on (hexagonal) polyominoes, by rotation or reflection. 
For any class $\F$ of polyominoes, with a monomial weigh $w$ 
corresponding to certain parameters, we  denote by
$|\F|_w$ the total weight (i.e. the  generating series) of this class.
If $\F$ is invariant under the action of $\D6$, then the set of orbits of this action
is denoted by $\F/\D6$. Burnside's Lemma enumerates these orbits in terms 
of the sets $\Fix(g)$ of fixed points of each of the elements $g$ of $\D6$, 
the \emph{symmetry classes} of $\F$. We write $\fix(g)=|\Fix(g)|_w$. 
Clearly we  have $\fix(\r)=\fix(\r^5),~\fix(\r^2)=\fix(\r^4)$ and, 
for symmetry reasons, $\fix(\da_1)=\fix(\da_2)=\fix(\da_3)$ and 
$\fix(\ds_1)=\fix(\ds_2)=\fix(\ds_3)$. In the following, we  will choose $\v=\da_2$, 
the vertical axis, and $\h=\ds_3$, the horizontal axis. 
We then have
\begin{eqnarray}
|\F/\D6|_w & = & {1\over12} \sum_{g\in\D6} \fix(g)    \nonumber \\
& = & {1\over12}\left(|\F|_w +2\,\fix(\r)+2\,\fix(\r^2)+\fix(\r^3)+3\,\fix(\v)+3\,\fix(\h)\right).
\label{eq:Burnside}
\end{eqnarray}
%
\subsubsection{The lattice of subgroups of $\D6$}
\begin{figure}[h]
\epsfysize=7cm
\centerline{\epsfbox{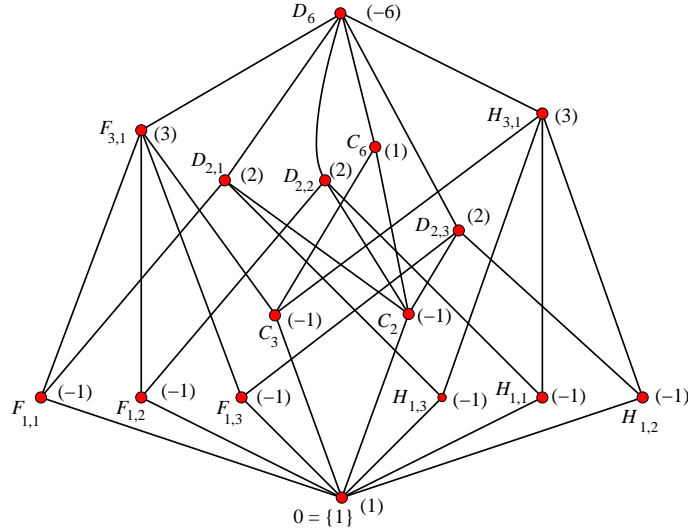}}
\vspace{-0.45\baselineskip}
\caption{The lattice of subgroups of $D_6$ ($\mu (0,H)$ in parenthesis)}
\label{fig:sgrD6}
\end{figure}

It is also possible to enumerate the convex polyominoes which are \emph{asymmetric} or which have
exactly the symmetries of a given subgroup $H$ of $\D6$, with the help of M\"obius 
inversion in the lattice of subgroups of $\D6$. This lattice and its  M\"obius function
are well described in Stockmeyer's Ph.D. thesis \cite{stock71}, for any dihedral group $\Dn$. 
We follow his nomenclature.
Apart from the trivial subgroups $0=\{\id\}$ and $1=\D6=<\r,\h>$, we  have the cyclic subgroups
$$ C_6=<\r>=\{1,\r ,\r ^2 ,\r ^3 ,\r ^4 ,\r ^5 \}, \quad
C_3=<\r^2>=\{1,\r ^2 ,\r ^4 \}\quad \mathrm{and} \quad
C_2=<\r^3>=\{1,\r ^3\}, $$
$$ F_{1,1}= <\ds_2>=\{1, \ds_2\}, \quad F_{1,2}= <\ds_1>=\{1, \ds_1\} \quad \mathrm{and} \quad 
F_{1,3}= <\h>=\{1, \h\}, $$
$$ H_{1,1}= <\da_3>=\{1, \da_3\},  \quad 
H_{1,2}= <\v>=\{1, \v\} \quad \mathrm{and} \quad H_{1,3}= <\da_1>=\{1, \da_1\}, $$
as well as the two generator subgroups 
$$ F_{3,1}= <\r^2, \ds_2>=\{1, \r ^2 , \r ^4 , \ds_1 , \ds_2 , \h\} = <\r^2, \h>, $$
$$ H_{3,1}= <\r^2, \da_3>=\{1, \r ^2 , \r ^4 , \da_1 , \v , \da_3\} = <\r^2, \v>, $$
and the $D_{2,j}= <\r^3, \tau\r^{2j}>$, $j=1,2,3$, that is
$$ D_{2,1}= <\r^3, \ds_2>=\{1, \r ^3, \ds_2, \da_1\},\quad 
\quad D_{2,2}= <\r^3, \ds_1>=\{1, \r ^3, \ds_1, \da_3\}, $$
$$\quad \quad \mathrm{and} \quad D_{2,3}= <\r^3, \h>=\{1, \r ^3, \h, \v\} = <\r^3, \h>. $$
The lattice of subgroups of $\D6$ is represented in Figure \ref{fig:sgrD6}.

For any subgroup $H$ of $\D6$ ($H\leq \D6$ ), we  set
\begin{equation} 
F_{\geq H}=\left| \{s\in \F~|~\mbox{stab}(s)\geq H\}\right|_{w} 
=\left| \{s\in \F~|~h\in H \Rightarrow h\cdot s=s\}\right|_{w}
\end{equation}
\noindent
and
\begin{equation} 
F_{=H}=\left| \{s\in \F~|~\mbox{stab}(s)= H\}\right|_{w} 
=\left| \{s\in \F~|~h\in H \Leftrightarrow h\cdot s=s\}\right|_{w}
\end{equation}

We have clearly, for any $H\leq \D6$
$$ F_{\geq H} = \sum_{H\leq K\leq \D6} F_{=K} $$
\noindent
and, by M\"obius inversion, 
$$ F_{=H}=\sum_{H\leq K\leq \D6} \mu(H,K) F_{\geq K} $$
In particular, 
the total weight of asymmetric polyominoes is given by 
\begin{equation}
F_{=0}=\sum_{K\leq \D6} \mu(0,K) F_{\geq K}
\label{eq:Fasym}
\end{equation}
Note that $F_{\geq 0}=|\F|_w$ and that for any cyclic subgroup $H=<h>$, $F_{\geq H}=\fix(h)$.
For reasons of symmetry (or by conjugation), we  have 
$F_{\geq D_{2,1}}= F_{\geq D_{2,2}} =  F_{\geq D_{2,3}}$. In the following, we  will take
$D_{2,3}=<\r^3, \h>$.  
The formula (\ref{eq:Fasym}) then yields
\begin{equation}
F_{=0}= |\F|_w -3\,\fix(\h) -3\,\fix(\v)-\fix(\r^2)-\fix(\r^3)+\fix(\r)
+6F_{\geq D_{2,3}}+3F_{\geq F_{3,1}}+3F_{\geq H_{3,1}}-6F_{\geq \D6}.
\label{eq:Fasymbis}
\end{equation}
For any subgroup $H$ of $\D6$, we  sometimes write 
$|\Fix(H)|_{q,t} = F_{\geq H}$ when the weight
is given by the area and the half perimeter. 
\vspace{.5cm}
\subsection{Growth phases of convex polyominoes}
\begin{figure}[h]
\epsfysize=5cm
\centerline{\epsfbox{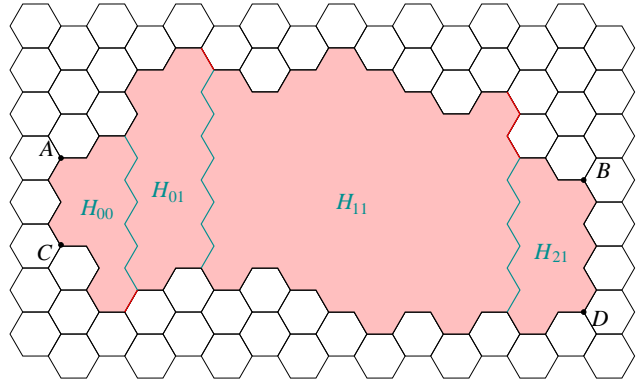}}
\vspace{-0.45\baselineskip}
\caption{Growth phases of a convex polyomino}
\label{fig:phasesconvex}
\end{figure}

Any convex polyomino can be decomposed into blocks according to the growth phases,
from left to right, of its upper and lower profiles. 
Figure \ref{fig:phasesconvex} gives an example of this decomposition. 
The upper profile is represented by the path from $A$ to $B$ along the  
upper boundary, and the lower profile, by the path from $C$ to $D$.
On the upper profile, we  speak of a \emph{weak growth} 
if the level rises by a half hexagon only with respect to the preceding column,
and of a \emph{strong} growth if the level rises by more than a half hexagon.
We define analogously a \emph{weak} or \emph{strong decrease}.
On the lower profile, a growth corresponds to a descent and a decrease,
to a rise.

The state in which a column lies is described by an ordered pair $(i,j)$, $i,j=0,1,2$;
the first component corresponds to the upper profile and the second, to the lower profile.   
The state $0$ corresponds to a (weak  or strong) growth, at the start of the polyomino,
the state $1$ to a weak growth or decrease, in an oscillation phase,
and the state $2$, to a strong or weak decrease, in the last part of the polyomino.
To pass from the state 0 to the state 1, there must be a first weak decrease, and to pass
from the state 0 or 1 to the state 2, there must be a strong decrease. Finally, the transitions 
from the state 1 to the state 0 and from the state 2 to the state 1 or 0 are impossible. 
Now, a block $H_{ij}$ is characterized by a maximal sequence of consecutive columns 
which are in the state $(i,j)$.

We can then view a convex polyomino as an assemblage of blocks and we will first enumerate
these blocks $H_{ij}$. In what follows, we give the various generating series
of the form $H_{ij}(x,q,u,v,t)$.
%
%
\subsubsection{The families $H_{00}$ and $H_{22}$}
The polyominoes of the classes $H_{00}$ and $H_{22}$ are easy to enumerate 
for they are in fact stack polyominoes. 
Here, only one of the two variables $u$ and $v$ is used at a time.
We have
\begin{equation}
H_{22}(x,q,u,t) = T(xt,ut^2,q)~~~\mathrm{and}~~~H_{00}(x,q,v,t) = T(xt,vt^2,q),
\label{eq:H00H22}
\end{equation}
where $T(x,u,q)$ is given by (\ref{eq:tas}).
%
%
\subsubsection{The families $H_{01}$, $H_{10}$, $H_{12}$, and $H_{21}$}
\begin{figure}[h]
\epsfysize=5cm
\centerline{\epsfbox{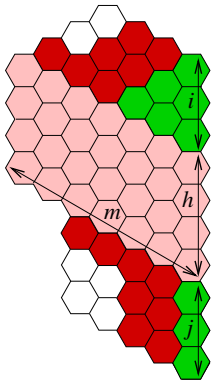}}
\vspace{-0.45\baselineskip}
\caption{Polyomino of $H_{10}$}
\label{fig:H10}
\end{figure}

The classes $H_{01}$, $H_{10}$, $H_{12}$, and $H_{21}$ of polyominoes are in 
bijection with each other by horizontal and vertical reflections and are thus equinumerous.
Figure \ref{fig:H10} shows a polyomino in $H_{10}$. We easily find that 
\begin{eqnarray} 
H_{10}(x,q,u,v,t)&=&\frac{xquvt^3}{1-quvt^2}+\frac{xt^2(1+qvt)}
{1-qvt^2}H_{10}(x,q,u,vq,t) \label{eq:H10} \\
&=&\sum_{m\geq 1}\frac{x^{m}q^{m}uvt^{2m+1}(-qvt;q)_{m-1}}
{(qvt^2;q)_{m-1}(1-q^{m}uvt^2)}.
\label{eq:H10bis}
\end{eqnarray}

Formula (\ref{eq:H10bis}) can be seen directly on Figure \ref{fig:H10}. We also see that
\begin{equation} 
H_{10}(x,q,u,v,t) = \sum_{h\geq 1}u^hv^h\sum_{m\geq 1}x^mq^{mh}
\sum_{i=0}^{m-1}v^iq^{{i+1}\choose 2}\left[ _{~~i}^{m-1} \right]_q
\sum_{j\geq 0}v^jq^jt^{2m+2h+i+2j-1}\left[ _{~~~j}^{m-2+j} \right]_q.
\label{eq:H10ter}
\end{equation}

Note that $H_{01}(x,q,u,v,t) = H_{10}(x,q,u,v,t)$ and that
$H_{12}(x,q,u,v,t) = H_{21}(x,q,u,v,t) = H_{10}(x,q,v,u,t)$.
%
%
\subsubsection{The families $H_{02}$ and $H_{20}$}
These two classes are in fact equivalent to staircase polyominoes:
\begin{equation}
H_{02}(x,q,u,v,t) = \Pa(x,q,u,v,t) = H_{20}(x,q,u,v,t).
\label{eq:H02H20}
\end{equation}
%
\subsubsection{The family $H_{11}$ \label{sec:H11}}
\begin{figure}[h]
\epsfysize=3.7cm
\centerline{\epsfbox{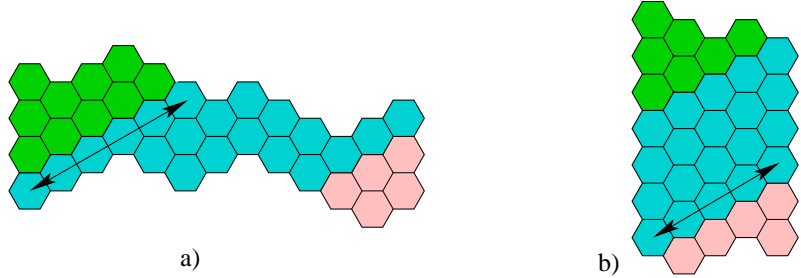}}
\vspace{-0.45\baselineskip}
\caption{Polyominoes of $H_{11}$}
\label{fig:H11}
\end{figure}

The class $H_{11}$ contains the convex polyominoes whose upper
and lower profiles are both oscillating. When we examine the diagonal row
of hexagons in the $\da_3$ axis (see Fig. \ref{fig:axesdereflexion}) 
containing the first column's lower cell, two subclasses of $H_{11}$ appear. In the first class,
denoted by $H_{11a}$, this diagonal row and those to its right (up to the last column) form
a staircase polyomino (rotated clockwise by a $\pi/3$ angle); see Figure \ref{fig:H11}a. 
In the second class, denoted by $H_{11b}$, this diagonal is the  basis of a rectangle
of height at least $2$; see Figure \ref{fig:H11}b.
In both cases, we find, above and below these objects (staircase
or rectangle), distinct part partitions which are left and right justified, respectively.

Recall that in the series $\Pa_{i,j}(x,q,t)$, defined by (\ref{eq:Gij}), the  variable $x$
marks the number of columns of the (unrotated) staircase polyomino. 
We rather use the link between
its width $\ell$, when rotated, and its half-perimeter $p$: $p=2\ell+1$.
Hence we have
\begin{equation} 
H_{11}(x,q,u,v,t) = H_{11a}(x,q,u,v,t) + H_{11b}(x,q,u,v,t)
\label{eq:H10AB}
\end{equation}
with
\begin{equation} 
H_{11a}(x,q,u,v,t) = \sum_{i\geq 1\,j\geq 1}x^{-{1\over2}}uv\,\Pa_{i,j}(1,q,tx^{1\over2})
D_{i-1}(ut,q)D_{j-1}(vt,q)
\label{eq:A(x,q,u,v,t)}
\end{equation}
and
\begin{equation} 
H_{11b}(x,q,u,v,t) = \sum_{n\geq 1}\frac{x^nq^{2n}u^2v^2t^{2n+3}D_{n-1}(ut,q)D_{n-1}(vt,q)}
{1-q^nuvt^2}.
\label{eq:B(x,q,u,v,t)}
\end{equation}
%
\section{Convex and directed convex polyominoes}
%
\subsection{Convex polyominoes \label{sec:convex}}
We denote by $C$, the class of all convex polyominoes and by $C_{ij}$, the subclass of polyominoes 
whose last column is in the state $(i,j),~i,j=0,1,2$. 
This determines a partition of $C$. 
To enumerate $C$, we must enumerate each of the classes $C_{ij}$. We give
the generating series $C_{ij}(x,q,v,t)$, 
using the growth phase decomposition of a convex polyomino, 
following essentially the method of Hassani \cite{ibn96}.

We use the notation $C_{ij}\otimes H_{i'j'}$ for the set of convex polyominoes
obtained by gluing together in all the possible ways a polyomino 
of $C_{ij}$ with one of $H_{i'j'}$.
We introduce the series $C_{ij,n}(x,q,t)$ and $H_{ij,n}(x,q,v,t)$ by the coefficient extractions
\begin{equation} 
C_{ij,n}(x,q,t) = [v^n]C_{ij}(x,q,v,t)~~\mathrm{and}~~H_{ij,n}(x,q,v,t) = [u^n]H_{ij}(x,q,u,v,t).
\label{eq:extractions}
\end{equation}
For example, we have $C_{00}=H_{00}$, $C_{10}=C_{00}\otimes H_{10}$ and
\begin{eqnarray} 
C_{10}(x,q,v,t) & = & \sum_{n\geq1}\left(\sum_{k=1}^n\frac{1}{t^{2k-1}}
C_{00,k}(x,q,t)\right)H_{10,n}(x,q,v,t)
\nonumber \\
& = & \sum_{n\geq1}\left(\sum_{k=1}^n tT_k(xt,q)\right)H_{10,n}(x,q,v,t)\\
& = & C_{01}(x,q,v,t).\nonumber
\label{eq:C10)}
\end{eqnarray}
Likewise, $C_{11}= (C_{00}+C_{10}+C_{01})\otimes H_{11} = C_{00}\otimes H_{11} +
C_{10}\otimes H_{11}  +C_{01}\otimes H_{11}$ and
\begin{equation} 
C_{00}\otimes H_{11}\,(x,q,v,t) = \sum_{n\geq 2}\frac{1}{t^{2n-2}}C_{00,n}(x,q,t)
H_{11,n-1}(x,q,v,t),
\label{eq:C00H11}
\end{equation}
\begin{equation} 
C_{10}\otimes H_{11}\,(x,q,v,t) 
= \sum_{n\geq 1}\frac{1}{t^{2n-1}}C_{10,n}(x,q,t)H_{11,n}(x,q,v,t) 
+ \sum_{n\geq 2}\frac{1}{t^{2n-2}}C_{10,n}(x,q,t)H_{11,n-1}(x,q,v,t).
\label{eq:C10H11}
\end{equation}

We have also $C_{02} = (C_{00}+C_{01})\otimes H_{02}$,
$$C_{12} = (C_{00}+C_{01}+C_{10}+C_{11}+C_{02})\otimes H_{12},$$
$$C_{22} = (C_{00}+C_{01}+C_{10}+C_{11}+C_{02}+C_{02})\otimes H_{22}.$$  
Finally,
\begin{equation} 
C(x,q,v,t) =(C_{00}+2C_{10}+C_{11}+2C_{02}+2C_{12}+C_{22})(x,q,v,t).
\label{eq:serieC}
\end{equation}
%
\subsection{Directed convex polyominoes
\label{sec:DC}}
\begin{figure}[h]
\epsfysize=5cm
\centerline{\epsfbox{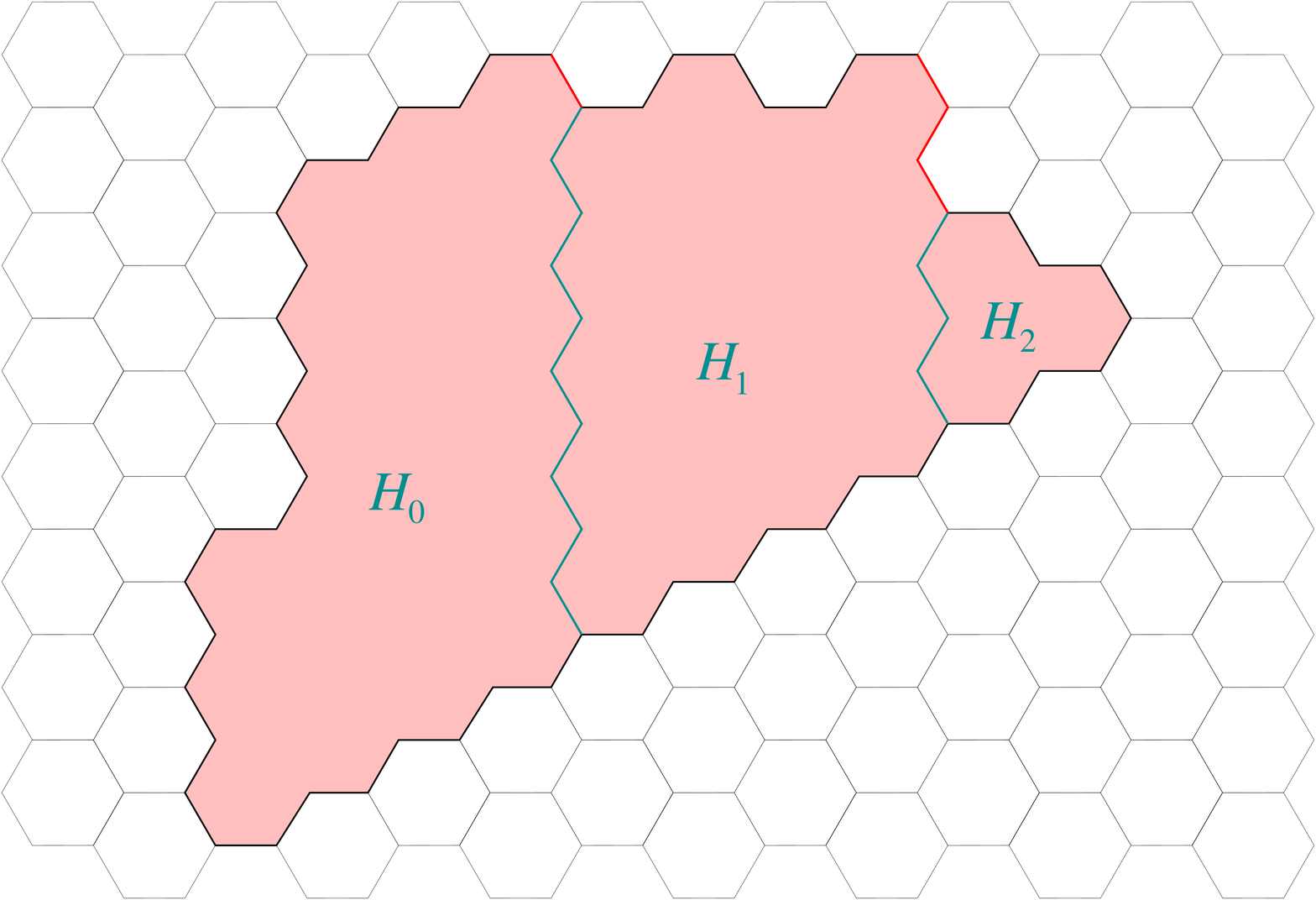}}
\vspace{-0.45\baselineskip}
\caption{Directed convex polyomino}
\label{fig:DC}
\end{figure}

There is a special class of convex polyominoes which will be particularly useful in the
following, namely those which are directed according to the North direction with 
a diagonal basis. See Figure \ref{fig:DC}. 
We call them simply \emph{directed convex}, and denote this class by $\DC$. 

\begin{figure}[h]
\epsfysize=5cm
\centerline{\epsfbox{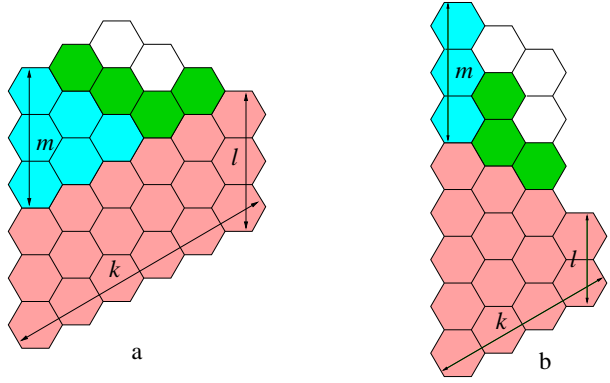}}
\vspace{-0.45\baselineskip}
\caption{Directed convex polyominos in $H_1$ and in $H_2$}
\label{fig:phasesDC}
\end{figure}

Analogously to convex polyominoes, any polyomino in $\DC$  can be 
decomposed into blocks $H_i$ according to the growth phases $i=0,1,~\mathrm{ou}~2$, 
of its upper profile. 
Their generating series can be computed directly by observation. 
Since the half perimeter can be deduced from the other parameters, 
the  variable $t$ does not appear here.  
For example, a polyomino in $H_0$ is identified with a partition and we find
\begin{equation}
H_{0}(x,q,v)
=\sum_{l\geq 1}\sum_{k\geq1}v^lx^kq^{l+k}\left[ _{~~l-1}^{l+k-2} \right]_q.
\label{eq:H0} 
\end{equation}
A polyomino in $H_1$ can be decomposed into pieces as shown in Figure \ref{fig:phasesDC}a, 
yielding 
\begin{equation}
H_1(x,q,u,v) = \sum_{l\geq 1}v^l\sum_{m\geq 0}u^{l+m}\sum_{k\geq m+1}
x^k q^{kl+\frac{(m)(m+1)}{2}}\left[ _{~m}^{k-1} \right]_q 
\label{eq:H1} .
\end{equation}
Likewise, for $H_2$, we find (see Figure \ref{fig:phasesDC}b)
\begin{equation}
H_2(x,q,u,v) = \sum_{l\geq 1}v^l\left( xq^lu^l + \sum_{k\geq2}x^k\sum_{m\geq0}u^{l+k+m-1} 
    q^{\frac{k(k+2l-1)}{2}+m} \left[ _{~~k-2}^{m+k-2} \right]_q \right).
\label{eq:H2}
\end{equation}

We denote by $\DC_i$, the class of directed convex polyominoes whose last column
is in the state $i, ~i=0,1,2$, and we introduce the notation
\begin{equation} 
\DC_{i,n}(x,q,t) = [v^n]\DC_{i}(x,q,v,t)~~\mathrm{and}~~H_{i,n}(x,q,v) = [u^n]H_{i}(x,q,u,v).
\label{eq:extractionsbis}
\end{equation}
We have 
\begin{equation}
\DC_{0}(x,q,v,t)={1\over t}H_{0}(xt^2,q,vt^2),
\label{eq:DC0} 
\end{equation}
$\DC_{1}=\DC_{0}\otimes H_{1}$ and
\begin{equation}
\DC_{1}(x,q,v,t)=\sum_{m\geq 1}{1\over t^m}\DC_{0,m+1}(x,q,t)H_{1,m}(xt^2,q,vt),
\label{eq:DC1} 
\end{equation}
and finally, $\DC_{2}=(\DC_{0}+\DC_{1})\otimes H_{2}$ and
\begin{equation}
\DC_{2}(x,q,v,t)=\sum_{m\geq 1}\sum_{h\geq2}{1\over t^m}
(\DC_{0,m+h}(x,q,t)+\DC_{1,m+h}(x,q,t))H_{2,m}(xt^2,q,vt).
\label{eq:DC2} 
\end{equation}
%

\section{Reflexive symmetry classes}
%
\subsection{Vertical symmetry \label{sec:vsym}}
\begin{figure}[h]
\epsfysize=4cm
\centerline{\epsfbox{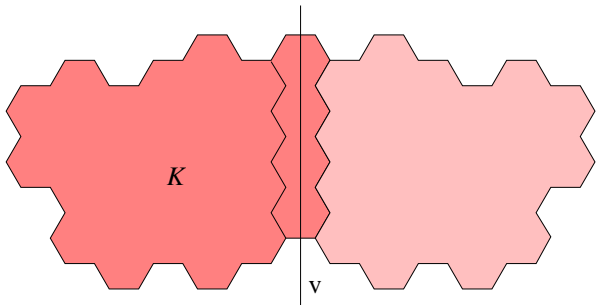}}
\vspace{-0.4\baselineskip}
\caption{Polyomino convex $v$-symmetric}
\label{fig:vsym}
\end{figure}

Consider a vertically symmetric ($\v$-symmetric) convex polyomino $P$. 
We see that the symmetry axis goes through a central column. 
Denote by $K$ the left fundamental region of $P$, including the central column.
See Figure \ref{fig:vsym}. We have 
\begin{eqnarray}
 K(x,q,v,t)&=& C_{00}(x,q,v,t)+2C_{10}(x,q,v,t)+C_{11}(x,q,v,t) \label{eq:serieK}\\
&=&\sum_{m\geq 1}K_m(x,q,t)v^m \nonumber
\end{eqnarray}
and
\begin{equation}
|\Fix(\v)|_{q,t}=\sum_{m\geq 1}\frac{1}{q^mt^{2m+1}}K_m(1,q^2,t^2).
\label{eq:fixv}
\end{equation}
%
\subsection{Horizontal symmetry \label{sec:hsym}}
\begin{figure}[h]
\epsfysize=6.5cm
\centerline{\epsfbox{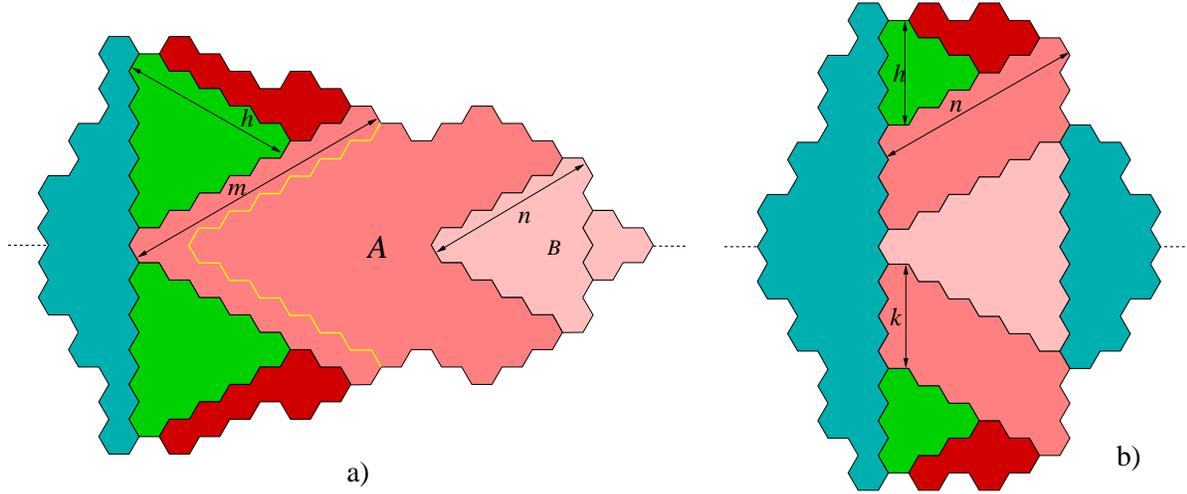}}
\vspace{-0.45\baselineskip}
\caption{$h$-symmetric convex polyominoes }
\label{fig:hsym}
\end{figure}

The class $S$ of $\h$-symmetric convex polyominoes is partitioned into three classes:
$S_a$ and $S_b$, wether or not we can find an \emph{arrowhead} polyomino
in the oscillating part (see Figures \ref{fig:hsym}a and \ref{fig:hsym}b) and 
the class $S_c$, if there does not exist an oscillating part. 

In order to construct a polyomino of the  class \emph{arrrowhead}, denoted by $A$,
we start with a triangle of side $n$, to which a symmetric stack is possibly
attached to form the $H_{22}$ phase; denote by $B$, this starting class of polyominoes.
From $B$, we construct $A$ by successively attaching $V$-shaped bands on the left, 
as illustrated in Figure \ref{fig:hsym}a. 
Let the variable $s$ mark the size of the last attached $V$'s upper left part.
We have
\begin{equation}
B(s,x,q,t)=sxqt^3+\sum_{n\geq 2}s^nx^nq^{n(n+1)/2}t^{3n} \TSO_{n-3}(xt,q)
\label{eq:B}
\end{equation}
and the generating series $A(s)=A(s,x,q,t)$ is 
characterized by the following functional equation, 
which can be solved by the usual method:
\begin{equation}
A(s) = B(s)+ \frac{s^2x^2q^3t^4}{1-sq^2} A(1)-A(sq^2).
\label{eq:AB}
\end{equation}
We set $A(s,x,q,t)=\sum_{m\geq 1}A_m(x,t,q)s^m$. 
To complete the polyomino, 
we must take into account the  parity of the  first oscillating column. 
The first case, illustrated in Figure \ref{fig:hsym}a,
is when this column is of odd size. In the second case, this even 
column is placed in front of the  arrrowhead. In conclusion, we obtain
\begin{eqnarray}
S_a(x,t,q) &=& \sum_{h\geq 0}q^{h(h+1)}t^{2h+2}\TS_{2h+2}(xt,q)
\sum_{m\geq h+1}\left[ _{~~h}^{m-1}\right]_{q^2}A_m(x,t,q) \nonumber \\
 & & \mbox{~~~~~~}+ ~\sum_{h\geq 1}q^{h(h+1)}t^{2h+3}\TS_{2h+1}(xt,q)
\sum_{m\geq h}\left[ _{\, h}^m\right]_{q^2}A_m(x,t,q).
\label{eq:Sa}
\end{eqnarray}

The computations for $S_b(x,t,q)$ and $S_c(x,t,q)$ are simpler. 
For $S_b$, there are also two parity cases and we find directly

$$
S_b(x,t,q) =
\sum_{n\geq 1}x^nq^{n+1\choose2}t^{3n}
\sum_{k\geq 1}q^{2kn}t^{4k}\TSO_{n+2k-3}(xt,q) 
\sum_{h=0}^{n-1}t^{2h+2}q^{h(h+1)}
\left[ _{~~h}^{n-1}\right]_{q^2}\TS_{2k+2h+2}(x,t,q) 
$$
\begin{equation}
\phantom{+}+~~\sum_{n\geq 0}x^nq^{n+1\choose2}t^{3n}
\sum_{k\geq 1}q^{2k(n+1)}t^{4k+1}\TSO_{n+2k-3}(xt,q) 
\sum_{h=0}^n t^{2h+2}q^{h(h+1)}\left[ _{\, h}^{n}\right]_{q^2}\TS_{2k+2h+1}(x,t,q)
\label{eq:Sb}
\end{equation}
and
\begin{equation}
S_c(x,t,q)=\sum_{h\geq 1}t^{2h}\TS_h(xt,q)\TSO_{h-3}(xt,q).
\label{eq:Sc}
\end{equation}
Finally,
\begin{equation}
|\Fix(\h)|_{q,t}=S_a(1,t,q)+S_b(1,t,q)+S_c(1,t,q).
\label{eq:fixh}
\end{equation}
%
\section{Rotational symmetry classes}
%
\subsection{Symmetry with respect to the $\pi/3$ radian rotation $\r$ \label{sec:rsym}}
The polyominoes which are symmetric with respect to the $\pi/3$ rotation ($\r$-symmetric)
are essentially formed of large hexagons decorated by stack polyominoes of the class $\TO$. We find
\begin{equation}
|\Fix(\r)|_{q,t}=\sum_{h\geq 1}t^{3(2h-1)}q^{3h(h-1)+1}\TO_{h-1}(t^6,q^6),
\label{eq:fixr}
\end{equation}
the series $\TO_{n}(x,q)$ being defined by (\ref{eq:recurrencandason}).
%
\subsection{Symmetry with respect to the $2\pi/3$ radian rotation $\r^2$ \label{sec:r2sym}}
The situation is more complex here. First we must distinguish the case where the rotation 
center is in the middle of an hexagon from the one where it is on a vertex.
This determines two subclasses, denoted by $\P$ and $\Q$.
%
\subsubsection{The rotation center is in the middle of an hexagon}
\begin{figure}[h]
\epsfysize=6cm
\centerline{\epsfbox{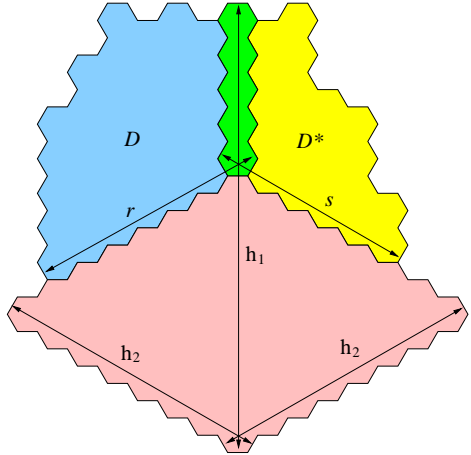}}
\vspace{-0.45\baselineskip}
\caption{Fundamental region of an $\r^2$-symmetric convex polyomino in $\P_1$}
\label{fig:r2sym}
\end{figure}

When the rotation center is in the middle of an hexagon, we consider the fundamental
region formed by the upper third of the $\r^2$-symmetric polyomino. The parameters
$h_1$ and $h_2$ are defined as the extent of the fundamental region in the
directions $\da_2 =\v$ and $\da_1$ (or $\da_3$), respectively.
There are three subcases according to wether $h_1>h_2$, $h_2>h_1$ or $h_1=h_2$, giving 
three subclasses denoted by $\P_1$, $\P_2$ and $\P_3$. Figure \ref{fig:r2sym} illustrates 
the first subcase $\P_1$.

In this figure, there is a basis formed of (one third of) a super-hexagon of 
\emph{radius} $h=h_2$, on top of which are placed a directed convex polyomino $D$
(see the subsection \ref{sec:DC}) and the image $D^{\ast}$ under $\v$ of another directed convex 
polyomino, these two polyominoes sharing a common column. 
Let $\ell$ be the  size of this common column, so that $h_1=h+\ell$. 
We must consider all such legal combinations $\DC \otimes \DC^{\ast}$
and take into account the added area and half perimeter over that of the super-hexagon.
For $i=0,1,2$, we set
$$\DC_{i,r,n}(q,t)={1\over t^{2r+n-1}}[x^r]\DC_{i,n}(x,q,t)={1\over t^{2r+n-1}}
   [x^r][v^n]\DC_{i}(x,q,v,t)$$
and $\DC_{r,n}(q,t) = (\DC_{0,r,n}+\DC_{1,r,n}+\DC_{2,r,n})(q,t)$.
We then have
\begin{eqnarray}
\P_1(q,t) & = & \sum_{h\geq 1}t^{3(2h-1)}q^{3h(h-1)+1} \nonumber \\
 & & \mbox{\hspace{0.5cm}}\sum_{l\geq 1}\sum_{r=1}^h\sum_{s=1}^h\frac{1}{q^{3l}}
\left(\DC_{0,r,l}(q^3,t^3)\DC_{0,s,l}(q^3,t^3)+2\DC_{0,r,l}(q^3,t^3)\DC_{1,s,l}(q^3,t^3)+ 
\right. \nonumber \\
 & & \left. \mbox{\hspace{3cm}}2\DC_{0,r,l}(q^3,t^3)\DC_{2,s,l}(q^3,t^3)+
\DC_{1,r,l}(q^3,t^3)\DC_{1,s,l}(q^3,t^3)\right).
\label{eq:P1qt}
\end{eqnarray}

For reasons of symmetry (a $\pi/3$ rotation), we see that $\P_2(q,t)=\P_1(q,t)$. 
\vspace{.05\baselineskip}

\begin{figure}[h]
\centerline{\input{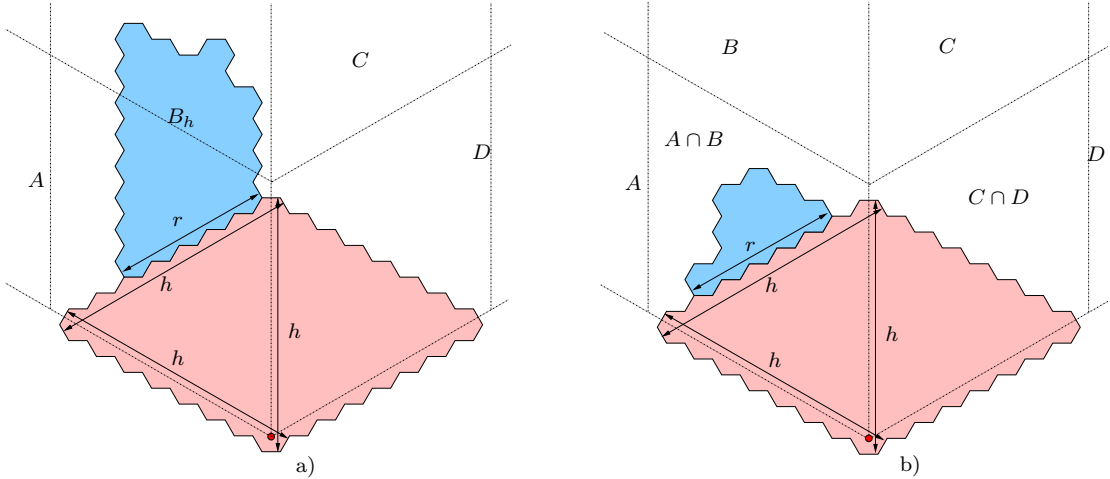}}
\vspace{-0.45\baselineskip}
\caption{Fundamental regions of $\r^2$-symmetric convex polyominoes in $\P_3$}
\label{fig:P3}
\end{figure}

Let us now consider the class $P_3$ of $\r^2$-symmetric convex polyominoes
with $h_1=h_2=h$. In this case, the added decorations on the super-hexagon can 
only occupy one of the sectors $A$, $B$, $C$ or $D$ shown in Figure \ref{fig:P3}, 
with the exception of the sectors $A\cap B$ and $C\cap D$ which can be simultaneously occupied.  

Denote by $B_h$ the class of admissible decorations in the sector $B$ over an hexagon
of \emph{side} $h$, and by $B_h(q,t)$ its generating series, where the variables $q$ and $t$
mark the added area and  half perimeter, respectively. We have
\begin{equation}
B_h(q,t)=\sum_{r=1}^{h-1}\left( (h-r)\left( q^{r}t+\DC_{1,r,1}(q,t)\right) +\DC_{2,r,1}(q,t)+
\sum_{j\geq 2}t^{j-1}\DC_{r,j}(q,t)\right)
\label{eq:Bhqt}
\end{equation}
The generating series will be  the same for the decorations located in the sectors
$A$, $C$ or $D$, for symmetry reasons. However, in the term $4B_h(q,t)$,
the decorations which, like the one shown in Figure \ref{fig:P3}b, are located in 
the intersection sectors $A\cap B$ or $C\cap D$, are counted twice. Observe 
that these decorations are in fact stack polyominoes of the class $\TO$ 
with generating series $\TO_{h-1}(t,q) -1$, whence the correcting term $-2(\TO_{h-1}(t,q) - 1)$.
Lastly, the term $(\TO_{h-1}(t,q) - 1)^2$ counts the simultaneous decorations in the sectors
$A\cap B$ and $C\cap D$  and the term 1 is added for the empty decoration.
Globally, we obtain
\begin{equation}
\P_3(q,t)=\sum_{h\geq 1}t^{3(2h-1)}q^{3h(h-1)+1}\left(4B_h(q^3,t^3)-
4TO_{h-1}(t^3,q^3)+TO^2_{h-1}(t^3,q^3)+4\right)
\label{eq:P3qt}
\end{equation}
%
\subsubsection{The rotation center is on a vertex}
\begin{figure}[h]
\epsfysize=5.5cm
\centerline{\epsfbox{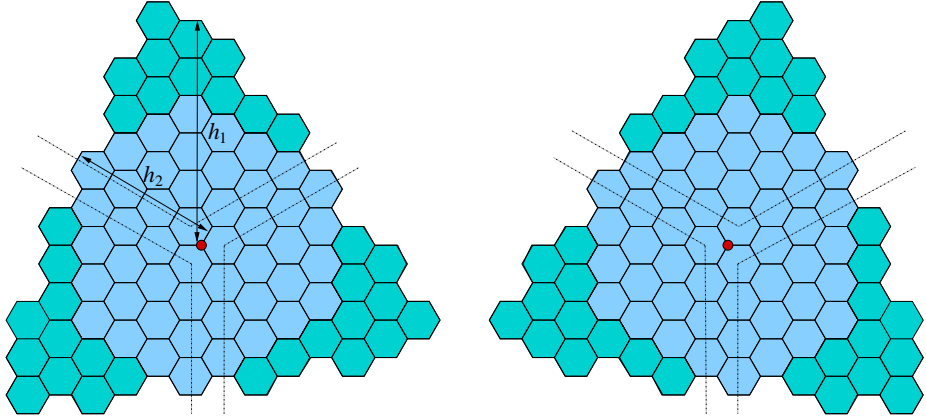}}
\vspace{-0.45\baselineskip}
\caption{$\r^2$-symmetric convex polyominoes in $\Q_1$}
\label{fig:Q1}
\end{figure}

The class of $\r^2$-symmetric convex polyominoes whose
rotation center is on a vertex is denoted by $\Q$.
There are three cases: $h_1>h_2$, $h_2>h_1$ and $h_1=h_2$ to which correspond
three series $\Q_1(q,t)$, $\Q_2(q,t)$ and $\Q_3(q,t)$ and two types
of central vertices as in Figure \ref{fig:Q1}. 
The computations are similar to the preceding case. The decorations are placed
over a \emph{pseudo-hexagon} and we find
\begin{eqnarray}
\Q_1(q,t) &=& 2\sum_{h\geq 1}t^{6h}q^{3h^2} 
\sum_{l\geq 1}\sum_{r=1}^{h+1}\sum_{s=1}^h\frac{1}
{q^{3l}}( \DC_{0,r,l}(q^3,t^3)\DC_{0,s,l}(q^3,t^3)+\DC_{0,r,l}(q^3,t^3)\DC_{1,s,l}
(q^3,t^3)+ \nonumber\\
& &\mbox{\hspace{2cm}}\DC_{0,r,l}(q^3,t^3)\DC_{2,s,l}(q^3,t^3)+
\DC_{1,r,l}(q^3,t^3)\DC_{1,s,l}(q^3,t^3)+ \nonumber\\
& & \mbox{\hspace{3cm}}\DC_{1,r,l}(q^3,t^3)\DC_{0,s,l}(q^3,t^3)+
\DC_{2,r,l}(q^3,t^3)\DC_{0,s,l}(q^3,t^3)) \label{eq:Q1qt}\\
&=&\Q_2(q,t) \nonumber
\end{eqnarray}
and
\begin{eqnarray}
\Q_3(q,t) &=& 2\sum_{h\geq 1}t^{6h}q^{3h^2}\left( 4+2B_h(q^3,t^3)+
2B_{h+1}(q^3,t^3)-2TO_{h-1}(t^3,q^3) \right. \nonumber\\
& & \mbox{\hspace{4cm}}\left. -2TO_{h}(t^3,q^3)+TO_{h-1}(t^3,q^3)TO_{h}(t^3,q^3)\right).
\label{eq:Q3qt}
\end{eqnarray}
%
\subsubsection{Global result}
Finally,
\begin{equation}
|\Fix(\r^2)|_{q,t} = 2\P_1(q,t)+\P_3(q,t)+2\Q_1(q,t)+\Q_3(q,t).
\label{eq:fixr2}
\end{equation}
%
%
\subsection{Symmetry with respect to the $\pi$ radian rotation $\r^3$}
\begin{figure}[h]
\epsfysize=5cm
\centerline{\epsfbox{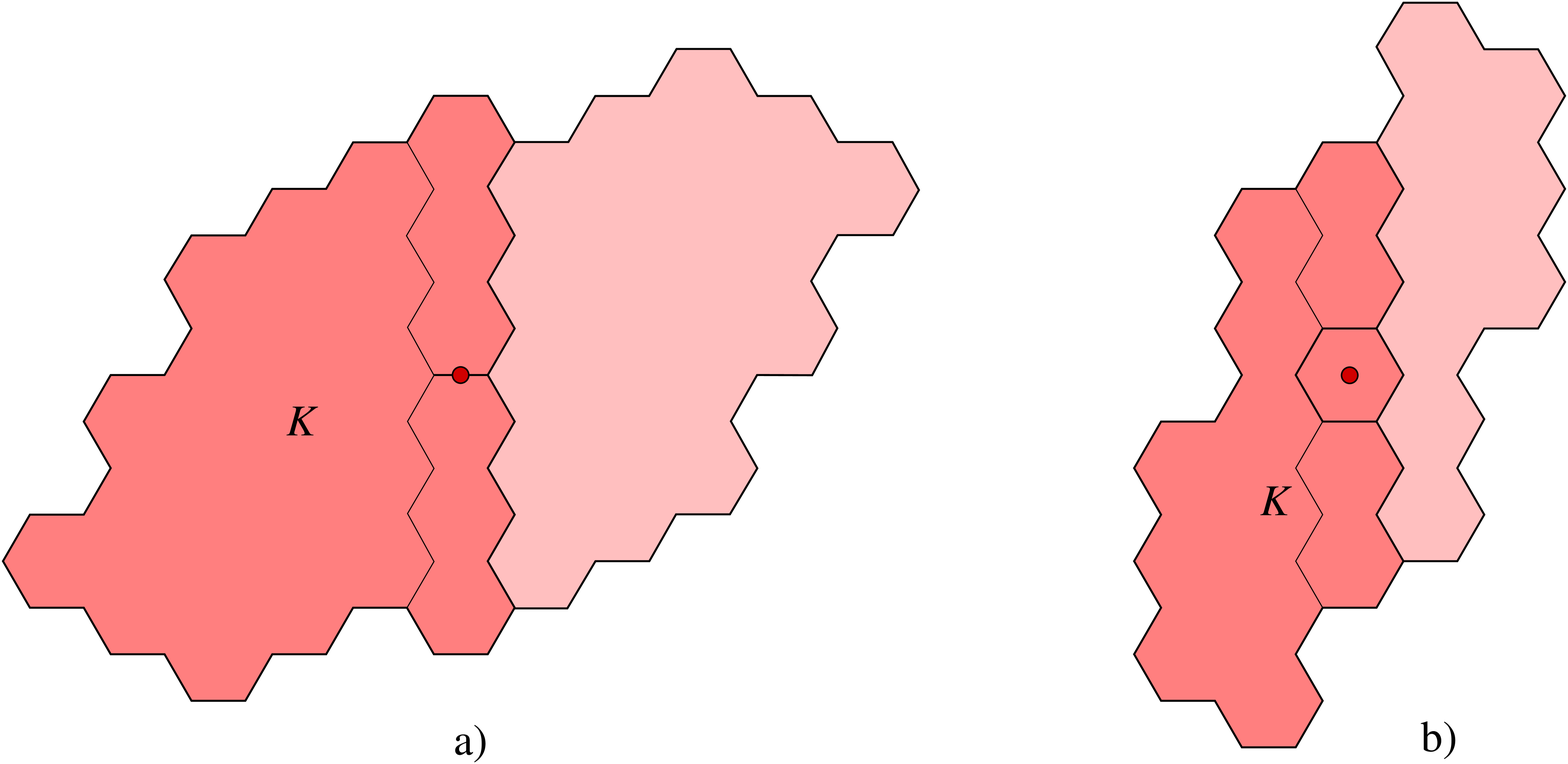}}
\vspace{-0.4\baselineskip}
\caption{$\r^3$-symmetric convex polyominos}
\label{fig:pisym}
\end{figure}

Here, the rotation center can be in the middle of an edge or of an hexagon. 
See Figure \ref{fig:pisym}.
If the rotation center is in the middle of an edge,
there are three similar cases corresponding to the three types of edges.
Consider the case of the horizontal edge and denote by $\CA$, the corresponding class.
Such a polyomino $P$ is shown in Figure \ref{fig:pisym}a.  
Denote by $K$ the left fundamental region of $P$, including the central column.
Note that this column is of even length.

If the rotation center is the middle of an hexagon, we denote by $\CH$ the corresponding class.
In this case, the central column is of odd length. See Figure \ref{fig:pisym}b.
The polyominoes $K$ which can occur as a fundamental region in one of these two cases are
\begin{equation}
K= C_{00} + C_{01} + C_{10} +C_{11} + C_{02} + C_{20}.
\label{eq:KA}
\end{equation}
Recall that the series $C_{ij,n}(x,q,t)$ is defined by equation (\ref{eq:extractions}),
with the index $n$ representing the size of the last column.
We then have
\begin{equation}
\CA(x,q,t)=\sum_{k\geq1}\frac{1}{xq^{2k}t^{4k+1}} 
(C_{00,2k} + 2C_{01,2k} +C_{11,2k} + 2C_{02,2k}) (x^2,q^2,t^2).
\label{eq:CAxqt}
\end{equation}
and 
\begin{equation}
\CH(x,q,t)=\sum_{k\geq 0}\frac{1}{xq^{2k+1}t^{4k+3}}
(C_{00,2k+1}+2C_{01,2k+1}+C_{11,2k+1}+2C_{02,2k+1})(x^2,q^2,t^2).
\label{eq:CHxqt}
\end{equation}
Finally,
\begin{equation}
|\Fix(\r^3)|_{q,t}=3\CA(1,q,t)+\CH(1,q,t).  \label{eq:fixr3qt}
\end{equation}
%
\section{Two generator symmetry classes}
%
\subsection{Symmetry with respect to $D_6$}

Since $D_6=<\r,\ds_2>$, convex polyominoes belonging to $\Fix(\D6)$ consist of 
super-hexagons with symmetric stack decorations (see the section \ref{sec:rsym}). 
We obtain
\begin{equation} 
|\Fix(D_6)|_{q,t}=\sum_{h\geq 1}t^{3(2h-1)}q^{3h(h-1)+1}\TSO_{h-1}(t^6,q^6).
\label{eq:fixD6qt}
\end{equation}
%
\subsection{Symmetry with respect to $F_{3,1}=<\r ^2,\ds_2>$}

We are guided by the $2\pi /3$ rotation symmetry class 
studied in section \ref{sec:r2sym}. The cases $h_1>h_2$ and  $h_2>h_1$ are impossible
because of the $\ds_2$-symmetry. There remains the case $h_1=h_2=h$. 
If the rotation center is in the middle of an hexagon (case $\P_3$), 
the sides of the superhexagon are decorated by symmetric stacks of type $\TSO$. 
Moreover the $\ds_2$-symmetry implies that the decorations 
are in the sectors $A\cap B$ and $C\cap D$.
By $\r^2$-symmetry, three of these stacks are identical and the three others
also, whence the formula $(TSO_{h-1}(t^3,q^3))^2$.
If the rotation center is a vertex (case $\Q_3$), we rather find
$TSO_{h-1}(t^3,q^3)TSO_{h}(t^3,q^3)$. Consequently,
\begin{eqnarray}
|\Fix(F_{3,1})|_{q,t}&=&\sum_{h\geq 1}t^{3(2h-1)}q^{3h(h-1)+1}(TSO_{h-1}(t^3,q^3))^2+\nonumber \\
& & \mbox{\hspace{4cm}} 2\sum_{h\geq 1}t^{6h}q^{3h^2}TSO_{h-1}(t^3,q^3)TSO_{h}(t^3,q^3).
\label{eq:fixr2dsqt}
\end{eqnarray}
%
\subsection{Symmetry with respect to $H_{3,1} = <\r ^2,\v>$}

We refer again to section \ref{sec:r2sym}.
The case where the center is a vertex is impossible because of the  vertical symmetry.
There remains the case where the center is an hexagon and the three subcases
$h_1>h_2$, $h_2>h_1$ and $h_1=h_2$ define three subclasses
$\R_1$, $\R_2$ and $\R_3$, respectively.
For the case where $h_1>h_2$, the  part $D^\ast$ of the  decoration (see Figure \ref{fig:r2sym})
is in fact the mirror image $\v\cdot D$ of $D$. 
Hence we obtain
\begin{eqnarray}
\R_1(q,t) &=& \sum_{h\geq 1}t^{3(2h-1)}q^{3h(h-1)+1}\sum_{l\geq 1}\sum_{r=1}^h\frac{1}{q^{3l}}
\left(C_{0,r,l}(q^6,t^6)+C_{1,r,l}(q^6,t^6)\right) \label{eq:R1qt}\\
&=&\R_2(q,t). \nonumber
\end{eqnarray}

If $h_1=h_2=h$, the decorations in the sectors $A\cap B$ and $C\cap D$
are mirror images of each other and we find
\begin{equation}
\R_3(q,t)=\sum_{h\geq 1}t^{3(2h-1)}q^{3h(h-1)+1}\TO_{h-1}(t^6,q^6).
\label{eq:R3qt}
\end{equation}
Finally,
\begin{equation}
|\Fix(H_{3,1})|_{q,t}=2\R_1(q,t)+\R_3(q,t).
\label{eq:FixH31qt}
\end{equation}
%
\subsection{Symmetry with respect to $D_{2,3} = <\r ^3,\h>$}
Observe that $D_{2,3}=<\h, \v>$. We thus refer to sections 
\ref{sec:vsym} and \ref{sec:hsym} on $\v$- and $\h$-symmetric polyominoes, respectively.
In order to obtain a $D_{2,3}$-symmetric convex polyomino, it suffices to take 
a $\v$-symmetric polyomino whose fundamental region $K$ (see the  figure \ref{fig:vsym})
is itself $\h$-symmetric.

The series $\CS_{00}(x,q,v,t)$, $\HS_{11}(x,q,u,v,t)$ and $\CS_{11}(x,q,v,t)$ 
are the $\h$-symmetric analogues of the series $C_{00}(x,q,v,t)$, $H_{11}(x,q,u,v,t)$ and 
$C_{11}(x,q,v,t)$ of sections \ref{sec:H11} and \ref{sec:convex}.
We have
\begin{equation}
\CS_{00,k}(x,q,t)=[v^k]\CS_{00}(x,q,v,t)=t^{2k}\TS_k(xt,q),
\label{eq:CS00xqvt}
\end{equation}
\begin{eqnarray}
\HS_{11}(x,q,u,v,t)&=&\frac{xquvt^3}{1-quvt^2}+xqvt^3~\HS_{11}(x,q,u,vq,t)+ \nonumber\\
& &\mbox{\hspace{3cm}}\frac{xt}{qv}\left( \HS_{11}(x,q,u,vq,t) 
-v(\frac{\HS_{11}(x,q,u,v,t)}{v})_{v=0}\right) \nonumber\\
&=& \sum_{k\geq 1}\HS_{11,k}(x,q,v,t)u^k, \label{eq:HS11xqvt}
\end{eqnarray}
\begin{eqnarray}
\CS_{11}(x,q,v,t)&=&\sum_{i\geq 2}
\frac{1}{t^{2i-2}}\CS_{00,i}(x,q,t)\HS_{11,i-1}(x,q,v,t) \label{eq:CS11xqvt}\\
&=& \sum_{k\geq 1}\CS_{11,k}(x,q,t)v^k, \nonumber
\end{eqnarray}
and finally
\begin{equation}
|\Fix(D_{2,3})|_{q,t} = \sum_{i\geq 1}
\frac{1}{q^it^{2i+1}}\left( \CS_{00,i}(1,q^2,t^2)+\CS_{11,i}(1,q^2,t^2)\right).
\end{equation}
%
\begin{table}[h]
\vspace{-0.25\baselineskip}
\caption{Symmetry classes of convex (hexagonal) polyominoes according to area
\label{table:area}}
\begin{center}
\begin{tabular}{|l|l|l|l|l|l|l|l|l|l|l|l|l|} \hline
Area & id & $h$ & $v$ & $r$ & $r^2$ & $r^3$ &  Orbits & $\D6$ & $F_{31}$ & 
$H_{31}$ & $D_{21}$ & Asym \\ \hline
1 & 1 & 1 & 1 & 1 & 1 & 1 & 1 & 1 & 1 & 1 & 1 & 0 \\ \hline
2 & 3 & 1 & 1 & 0 & 0 & 3 & 1 & 0 & 0 & 0 & 1 & 0 \\ \hline
3 & 11 & 3 & 3 & 0 & 2 & 3 & 3 & 0 & 2 & 0 & 1 & 0 \\ \hline 
4 & 38 & 2 & 4 & 0 & 2 & 12 & 6 & 0 & 0 & 2 & 2 & 24 \\ \hline
5 & 120 & 6 & 10 & 0 & 0 & 12 & 15 & 0 & 0 & 0 & 2 & 72 \\ \hline
6 & 348 & 6 & 12 & 0 & 6 & 42 & 38 & 0 & 2 & 0 & 2 & 264 \\ \hline
7 & 939 & 9 & 27 & 1 & 3 & 37 & 91 & 1 & 1 & 3 & 3 & 816 \\ \hline
8 & 2412 & 12 & 30 & 0 & 0 & 126 & 222 & 0 & 0 & 0 & 4 & 2184 \\ \hline
9 & 5973 & 17 & 63 & 0 & 12 & 99 & 528 & 0 & 0 & 0 & 3 & 5640 \\ \hline
10 & 14394 & 20 & 66 & 0 & 6 & 336 & 1250 & 0 & 2 & 4 & 4 & 13836 \\ \hline
11 & 34056 & 30 & 142 & 0 & 0 & 252 & 2902 & 0 & 0 & 0 & 6 & 33324 \\ \hline
12 & 79602 & 38 & 140 & 0 & 18 & 840 & 6751 & 0 & 2 & 0 & 4 & 78240 \\ \hline
13 & 184588 & 46 & 310 & 1 & 13 & 616 & 15525 & 1 & 1 & 5 & 8 & 182952 \\ \hline
14 & 426036 & 62 & 286 & 0 & 0 & 2028 & 35759 & 0 & 0 & 0 & 8 & 423012 \\ \hline
15 & 980961 & 69 & 665 & 0 & 30 & 1461 & 82057 & 0 & 2 & 0 & 7 & 977316 \\ \hline
16 & 2256420 & 100 & 580 & 0 & 18 & 4788 & 188607 & 0 & 0 & 6 & 8 & 2249640 \\ \hline
17 & 5189577 & 115 & 1441 & 0 & 0 & 3435 & 433140 & 0 & 0 & 0 & 11 & 5181540 \\ \hline
18 & 11939804 & 154 & 1184 & 0 & 50 & 11142 & 996255 & 0 & 2 & 0 & 12 & 11924676 \\ \hline
19 & 27485271 & 175 & 3145 & 1 & 27 & 8005 & 2291941 & 1 & 1 & 7 & 13 & 27467376 \\ \hline
20 & 63308532 & 238 & 2458 & 0 & 0 & 25800 & 5278535 & 0 & 0 & 0 & 16 & 63274740 \\ \hline
\end{tabular}
\end{center}
\vspace{-0.95\baselineskip}
\end{table}
%
%
\section{Conclusion}
It is now possible to use Burnside's formula (\ref{eq:Burnside}), with $\F=C$,
to enumerate the free convex polyominoes, 
according to area and half perimeter. Some numerical results are given in
tables 1 and 2, according to area only (up to area 20) or to half perimeter only
(up to half perimeter 16). See under the column "Orbits".

It is also possible to enumerate asymmetric convex polyominoes 
with the help of formula (\ref{eq:Fasymbis}), with $\F=C$. 
Some results are found in the tables 1 and 2.
It is clear on these tables that almost all convex polyominoes are asymmetric.

All these numerical results were verified experimentally by an 
exhaustive computerized enumeration.
%
\begin{table}[h]
\vspace{-0.25\baselineskip}
\caption{Symmetry classes of convex (hexagonal) polyominoes according to half perimeter
\label{table:demiperi}}
\begin{center}
\begin{tabular}{|l|l|l|l|l|l|l|l|l|l|l|l|l|} \hline
${1\over2}$per. & id & $\h$ & $\v$ & $\r$ & $\r^2$ & $\r^3$ &  Orbits & $\D6$ & $F_{31}$ & 
$H_{31}$ & $D_{21}$ & Asym \\ \hline
3 & 1 & 1 & 1 & 1 & 1 & 1 & 1 & 1 & 1 & 1 & 1 & 0 \\ \hline
4 & 0 & 0 & 0 & 0 & 0 & 0 & 0 & 0 & 0 & 0 & 0 & 0 \\ \hline
5 & 3 & 1 & 1 & 0 & 0 & 3 & 1 & 0 & 0 & 0 & 1 & 0 \\ \hline 
6 & 2 & 2 & 0 & 0 & 2 & 0 & 1 & 0 & 2 & 0 & 0 & 0 \\ \hline
7 & 12 & 2 & 4 & 0 & 0 & 6 & 3 & 0 & 0 & 0 & 2 & 0 \\ \hline
8 & 18 & 2 & 0 & 0 & 0 & 0 & 2 & 0 & 0 & 0 & 0 & 12 \\ \hline
9 & 59 & 5 & 9 & 1 & 5 & 19 & 11 & 1 & 3 & 3 & 3 & 24 \\ \hline
10 & 120 & 8 & 0 & 0 & 0 & 0 & 12 & 0 & 0 & 0 & 0 & 96 \\ \hline
11 & 318 & 10 & 24 & 0 & 0 & 48 & 39 & 0 & 0 & 0 & 6 & 204 \\ \hline
12 & 714 & 14 & 0 & 0 & 12 & 0 & 65 & 0 & 4 & 0 & 0 & 672 \\ \hline
13 & 1743 & 25 & 59 & 0 & 0 & 129 & 177 & 0 & 0 & 0 & 7 & 1368 \\ \hline
14 & 4008 & 36 & 0 & 0 & 0 & 0 & 343 & 0 & 0 & 0 & 0 & 3900 \\ \hline
15 & 9433 & 53 & 143 & 2 & 28 & 323 & 867 & 2 & 6 & 8 & 15 & 8616 \\ \hline
16 & 21672 & 76 & 0 & 0 & 0 & 0 & 1825 & 0 & 0 & 0 & 0 & 21444 \\ \hline
\end{tabular}
\end{center}
\vspace{-0.75\baselineskip}
\end{table}
%
{\small%
\bibliographystyle{plain}
\bibliography{hexagonal}}
%
\end{document}